\input amstex\documentstyle {amsppt}  
\pagewidth{12.5 cm}\pageheight{19 cm}\magnification\magstep1
\topmatter
\title Convolution of almost characters\endtitle
\author G. Lusztig\endauthor
\address Department of Mathematics, M.I.T., Cambridge, MA 02139\endaddress
\thanks Supported in part by the National Science Foundation\endthanks
\endtopmatter   
\document

\define\op{\oplus}

\define\m{\mapsto}

\define\bxt{\boxtimes}
\define\T{\times}

\define\nl{\newline}
\redefine\i{^{-1}}
\define\fra{\frac}
\define\un{\underline}
\define\ov{\overline}
\define\ot{\otimes}

\define\Hom{\text{\rm Hom}}

\define\tr{\text{\rm tr}}

\define\g{\gamma}

\define\e{\epsilon}
\define\et{\eta}

\define\p{\pi}

\define\r{\rho}
\define\s{\sigma}

\redefine\G{\Gamma}
\redefine\D{\Delta}

\define\qq{\bold q}

\define\CC{\bold C}

\define\FF{\bold F}

\define\NN{\bold N}

\define\ZZ{\bold Z}

\define\cc{\Cal C}

\define\cg{\Cal G}

\define\cm{\Cal M}

\define\cu{\Cal U}

\define\uc{\un c}
\define\Ve{\text{\rm Vec}}
\define\hs{\hat\s}
\define\cguc{\cg_{\uc}}
\define\ccuc{\cc_{\uc}}
\define\UR{L1}
\define\CR{L2}
\define\CS{L3}
\define\IC{L4}
\define\LC{L5}
\subhead 0\endsubhead
In this paper we show that the convolution product of  "almost characters" of a 
connected reductive group over a finite field is given by "structure constants" whose 
leading coefficients can be interpreted in $K$-theoretic terms and in particular are
natural numbers. A conjectural analogue for character sheaves is also stated.

\subhead 1\endsubhead
Let $\cg$ be a finite group. Let $\Ve_\cg$ be the category whose objects are 
$\CC$-vector bundles on $\cg$ which are $\cg$-equivariant for the conjugation action of
$\cg$ on $\cg$. For $U\in\Ve_\cg$ let $U_g$ denote the fibre of $U$ at $g\in\cg$; the 
$\cg$-equivariant structure on $U$ is given by isomorphisms $\psi^U_x:U_y@>>>U_{xyx\i}$
for any $x,y\in\cg$. For $U,U'\in\Ve_\cg$ we define $U*U'\in\Ve_\cg$ as follows: for
$y\in\cg$, we have $(U*U')_y=\op_{(g,g')\in\cg\T\cg;gg'=y}U_g\ot U'_{g'}$; the 
$\cg$-equivariant structure on $U*U'$ is defined by 
$$\psi^{U*U'}_x:\op_{(g,g')\in\cg\T\cg;gg'=y}U_g\ot U'_{g'}@>>>
\op_{(h,h')\in\cg\T\cg;hh'=xyx\i}U_h\ot U'_{h'}$$
given by the direct sum over $\{(g,g')\in\cg\T\cg;gg'=y\}$ of
$\psi^U_x\ot\psi^{U'}_x:U_g\ot U'_{g'}@>>>U_{xgx\i}\ot U'_{xg'x\i}$. Then $\Ve_\cg$ 
becomes a tensor category. In \cite{\LC, p.242} it is shown that $U*U'$ is isomorphic 
to $U'*U$ for any $U,U'\in\Ve_\cg$. More directly, for any $y\in\cg$ we define an
isomorphism 
$$\p_y:\op_{(g,g')\in\cg\T\cg;gg'=y}U_g\ot U'_{g'}@>>>
\op_{(h,h')\in\cg\T\cg;hh'=y}U'_h\ot U_{h'}$$
by $\p_y(u\ot u')=\psi^{U'}_g(u')\ot u\in U'_{gg'g\i}\ot U_g$ (where
$u\in U_g,u'\in U'_{g'}$). Then $\p=(\p_y)_{y\in\cg}:U*U'@>>>U'*U$ is an isomorphism of
$\cg$-equivariant vector bundles.

Let $\cm(\cg)$ be the set of pairs $(g,\r)$ where $g\in\cg$ is defined up to
$\cg$-conjugacy and $\r$ is an irreducible representation of the centralizer $Z_\cg(g)$
defined up to isomorphism. As in \cite{\UR} we have a pairing 
$\{,\}:\cm(\cg)\T\cm(\cg)@>>>\CC$ given by
$$\{(g,\r),(g',\r')\}=|Z_\cg(g)|\i|Z_\cg(g')|\i\sum_{h\in\cg;h\i gh\in Z_\cg(g')}
\tr(hg'h\i,\r)\tr(h\i g\i h,\r').$$
An object in $\Ve_\cg$ is said to be irreducible if it is non-zero and is not a direct
sum of two non-zero objects of $\Ve_\cg$. According to \cite{\LC,2.5}, the set of
isomorphism classes of irreducible objects of $\Ve_\cg$ is in natural bijection with
the set $\cm(\cg)$: to $m=(g,\r)\in\cm(\cg)$ corresponds an object $U=U^m\in\Ve_\cg$ 
such that $U_h$ is zero for $h\in\cg$ not conjugate to $g$ and 
$\tr(\psi^U_x:U_g@>>>U_g)=\tr(x,\r)$ for any $x\in Z_\cg(g)$. Then 
$(U^m)_{m\in\cm(\cg)}$ is a $\ZZ$-basis for the Grothendieck group $K_\cg(\cg)$ of 
$\Ve_\cg$. For $m,m'$ in $\cm(\cg)$, we can write
$$U^m*U^{m'}\cong\op_{m''\in\cm(\cg)}(U^{m''})^{\op p_{m,m',m''}}$$
where $p_{m,m',m''}\in\NN$ are multiplicities. Now $*$ induces a structure of 
commutative ring on $K_\cg(\cg)$ in which $p_{m,m',m''}$ are structure constants:
$$U_mU_{m'}=\sum_{m''\in\cm(\cg)}p_{m,m',m''}U_{m''}.\tag a$$
According to \cite{\LC, 2.5(d)} for any $(z,\et)\in\cm(\cg)$, the assignment
$$U^m\m(\dim\et)\i|Z_\cg(z)|\{m,(z,\et)\}$$
defines a ring homomorphism $K_\cg(\cg)@>>>\CC$. Applying this ring homomorphism to (a)
gives
$$(\dim\et)\i|Z_\cg(z)|\{m,(z,\et)\}\{m',(z,\et)\}=\sum_{m''\in\cm(\cg)}
p_{m,m',m''}\{m'',(z,\et)\}.$$
Since the matrix $(\{m,m'\})$ has square $1$ (see \cite{\UR}) we deduce that
$$p_{m,m,m''}=\sum_{(z,\et)\in\cm(\cg)}(\dim\et)\i|Z_\cg(z)|\{m,(z,\et)\}\{m',(z,\et)\}
                                                               \{(z,\et),m''\}.\tag b$$
\subhead 2\endsubhead
Let $\G$ be a finite group and let $f,f'$ be two class functions $\G@>>>\CC$. Recall 
that the convolution $f*f'$ is the class function $\G@>>>\CC$ given by
$(f*f')(\g)=\sum_{\g_1\in\G}f(\g_1)f'(\g_1\i\g)$. It is well known that given 
irreducible characters $f,f'$ of $\G$, we have $f*f'=0$ if $f\ne f'$ and 
$f*f=(|\G|/f(1))f$.

\subhead 3\endsubhead
Let $G$ be a connected reductive algebraic group over an algebraic closure of a finite 
field $\FF_q$, with a fixed $\FF_q$-rational structure. The vector space of class
functions $G(\FF_q)@>>>\CC$ has an orthonormal basis consisting of the irreducible 
characters and also an orthonormal basis consisting of "almost characters" (see 
\cite{\CR, 13.6}), closely related to character sheaves. We would like to study the 
convolution of two almost characters. Since the basis of almost characters differs from
the basis of irreducible characters by an explicitly known almost diagonal matrix, we 
see that the convolution of two almost characters is a linear combination of a small
number of almost characters. We want to make explicit the coefficients in this linear
combination. For simplicity, here we restrict ourselves to "unipotent almost 
characters"; these are the almost characters that are linear combinations of unipotent 
characters (they form a basis of the space spanned by the unipotent characters). We
also assume that $G$ is split over $\FF_q$ and that $G$ modulo its centre is simple.

To each two sided cell $\uc$ of the Weyl group $W$ of $G$ we attach as in 
\cite{\CR, \S4} a certain finite group $\cguc$. As in \cite{\CR} the set of unipotent 
characters (resp. unipotent almost characters) of $G^F$ can be indexed as 
$(\s_{\uc,m})$ (resp. $(\hs_{\uc,m})$) where $\uc$ is a two sided cell of $W$, 
$m\in\cm(\cguc)$ and
$$\hs_{\uc,m}=\sum_{m'\in\cm(\cguc)}\{m',m\}\D(\uc,m')\s_{\uc,m'}\tag a$$
$$\s_{\uc,m}=\sum_{m'\in\cm(\cguc)}\{m',m\}\D(\uc,m)\hs_{\uc,m'}\tag b$$
with $\D(\uc,m')=\pm 1$. Moreover, to any two-sided cell $\uc$ one can attach an 
integer $A_{\uc}\ge 0$ such that for any $m=(g,\r)\in\cguc$, $|G(\FF_q)|/\s_{\uc,m}(1)$
is a polynomial in $q$ with constant rational coefficients of the form
$$\fra{|Z_{\cguc}(g)|}{\dim\r}q^{D-A_{\uc}}+\text{ lower powers of }q.\tag c$$
(Here $D=\dim G$.)

We now consider two unipotent almost characters $\hs_{\uc,m},\hs_{\uc',m'}$. From (a)
we see that 

$\hs_{\uc,m}*\hs_{\uc',m'}=0$ if $\uc\ne\uc'$.
\nl
Now assume that $\uc=\uc'$. From (a),(b) we have
$$\align&\hs_{\uc,m}*\hs_{\uc,m'}\\&=\sum_{m_1,m_2\in\cm(\cguc)}\{m_1,m\}\{m_2,m'\}
\D(\uc,m_1)\D(\uc,m_2)\s_{\uc,m_1}*\s_{\uc,m_2}\\&
=\sum_{m_1\in\cm(\cguc)}\{m_1,m\}\{m_1,m'\}|G(F_\qq)|\s_{\uc,m_1}(1)\i\s_{\uc,m_1}\\&
=\sum_{m_1\in\cm(\cguc)}\{m_1,m\}\{m_1,m'\}|G(F_\qq)|\s_{\uc,m_1}(1)\i\\&
\sum_{m''\in\cm(\cguc)}\{m'',m_1\}\D(\uc,m_1)\hs_{\uc,m''}.\endalign$$
Thus
$$\hs_{\uc,m}*\hs_{\uc,m'}=\sum_{m'_1\in\cm(\cguc)}N_{m,m',m''}\hs_{\uc,m''}$$ 
where, for $m,m',m''\in\cm(\cguc)$ we have
$$N_{m,m',m''}=\sum_{m_1\in\cm(\cguc)}\{m_1,m\}\{m_1,m'\}\{m'',m_1\}\D(\uc,m_1)
|G(F_\qq)|\s_{\uc,m_1}(1)\i.$$
Using (c) we see that $N_{m,m',m''}$ is a polynomial in $q$ with constant rational 
coefficients of the form
$$n_{m,m',m''}q^{D-A_{\uc}}+\text{ lower powers of }q$$
where
$$n_{m,m',m''}=\sum_{(z,\et)\in\cm(\cguc)}\fra{|Z_{\cguc}(z)|}{\dim\et}
\D(\uc,(z,\et))\{(z,\et),m\}\{(z,\et),m'\}\{m'',(z,\et)\}.$$
If $\D(\uc,m_1)$ is $1$ for any $m_1\in\cm(\cguc)$ then, using 1(b) and the identity 
$\{m_1,m_2\}=\ov{\{m_2,m_1\}}$, we see that 
$$n_{m,m',m''}=p_{m,m',m''}$$
where $p_{m,m',m''}$ is defined as in Section 1 in terms of $\cg=\cguc$. (We use that
the integer $p_{m,m',m''}$ is fixed by complex conjugation.) Thus, the leading 
coefficient $n_{m,m',m''}$ has a $K$-theoretic interpretation; in particular, it is in 
$\NN$.

A similar result holds when $m_1\m\D(\uc,m_1)$ is not identically $1$. (Then $G$ is of 
type $E_7$ and $\uc$ is a two-sided cell with $|\uc|=2\T 512^2$ or $G$ is of type $E_8$
and $\uc$ is a two-sided cell with $|\uc|=2\T 4096^2$.) In this case $\cguc=\ZZ/2\ZZ$, 
$\cm_{\cguc}$ may be identified with the (additive) abelian group 
$\ZZ/2\ZZ\T\Hom(\ZZ/2\ZZ,\CC^*)$ and $\D(\uc,m_1)=2\{m_0,m_1\}$ where $m_0=(0,\e)$ and 
$\e\ne 0$. Moreover, in this case we have 
$\{m_1,m_2\}\{m_1,m_3\}=\{m_1,m_2+m_3\}/2$ and $\{m_1,m_2\}=\{m_2,m_1\}$. It follows 
that
$$n_{m,m',m''}=\sum_{m_1\in\cm(\cguc)}\{m_1,m+m'+m''+m_0\}/2.$$
On the other hand we have 
$$p_{m,m',m''}=\sum_{m_1\in\cm(\cguc)}\{m_1,m+m'+m''\}/2.$$
Hence 
$$n_{m,m',m''}=p_{m+m_0,m',m''}=p_{m,m'+m_0,m''}=p_{m,m',m''+m_0}.$$
{\it Remark.} If $G$ is a classical group, then for any $\uc$, $\cm(\cg_\cu)$ may be
regarded as an $\FF_2$-vector space of even dimension and for $m,m',m''$ in
$\cm(\cg_\cu)$ we have $n_{m,m',m''}=p_{m,m',m''}=1$ if $m+m'+m''=0$ and
$n_{m,m',m''}=p_{m,m',m''}=0$ if $m+m'+m''\ne 0$.

\subhead 4\endsubhead
Under some mild restriction on the characteristic of $\FF_q$, the "unipotent" character
sheaves on $G$ can again be indexed as $(K_{\uc,m})$) where $\uc$ is a two sided cell 
of $W$ and $m\in\cm(\cguc)$ (see \cite{\CS}). For $\uc$ as above, let $\ccuc$ be the 
category whose objects are perverse sheaves which are direct sums of character sheaves 
of the form $K_{\uc,m}$ for various $m\in\cm(\cguc)$. Let $\mu:G\T G@>>>G$ be 
multiplication. For $K,K'\in\ccuc$ we set
$K\circ K'={}^pH^{2D-2A_{\uc}}(\mu_!(K\bxt K'))$ (notation of Section 3; ${}^pH^i()$
denotes the $i$-th perverse cohomology sheaf).

We conjecture that $K\circ K'\in\ccuc$ and that $K,K'\m K\circ K'$ defines on $\ccuc$ a
structure of tensor category; moreover, this tensor category should be equivalent to 
$\Ve_{\cguc}$. In particular, for $m,m'\in\cm(\cguc)$, we should have
$$K_{\uc,m}\circ K_{\uc,m'}\cong\op_{m''\in\cm(\cguc)}K_{\uc,m''}^{\op p_{m,m',m''}}$$
with $p_{m,m',m''}$ defined as in Section 1 in terms of $\cg=\cguc$.

This conjecture is suggested by the results in Section 3. We expect that it extends 
also to the non-unipotent case and that its proof should involve the geometric 
interpretation of the indexing of the character sheaves in $\ccuc$ by irreducible 
vector bundles in $\Ve_{\cguc}$ indicated in \cite{\IC, 4.7(a)}.

\enddocument